\documentclass[12pt]{amsart}
\usepackage{amscd}     
\usepackage{amssymb}
\usepackage{amsmath, amsthm, graphics, dsfont}
\usepackage{xypic}     
\LaTeXdiagrams        
\usepackage[all]{xy}
\xyoption{2cell} \UseAllTwocells \xyoption{frame} \CompileMatrices
\allowdisplaybreaks[3]
\usepackage{amsfonts}
\usepackage[normalem]{ulem}
\usepackage{hyperref}
\usepackage{tikz}
\usepackage{mathtools}
\usepackage{verbatim} 
\usepackage{MnSymbol}

\usepackage{latexsym}
\usepackage{epsfig}
\usepackage{enumerate}
\usepackage{times}

\newtheorem{prop}{Proposition}

\newtheorem{thm}[prop]{Theorem}
\newtheorem{rem}[prop]{Remark}

\newtheorem{rmk}{Remark}

\newcommand{\Hilb}{\mathsf{Hilb}^n(S)}
\newcommand{\Sym}{\mathsf{Sym}^n(S)}

\newcommand{\Mbar}{\overline{\mathcal{M}}}
\newcommand{\ev}{\text{ev}}
\newcommand{\K}{\mathcal{K}}

\title{On $\mathsf{Hilb}/\mathsf{Sym}$ correspondence for K3 surfaces}

\author[Genlik]{Deniz Genlik}
\address{Department of Mathematics\\University of Illinois at Urbana-Champaign\\1409 W. Green Street (MC-382)\\Urbana, IL 61801\\ USA}
\email{genlik@illinois.edu}

\author[Tseng]{Hsian-Hua Tseng}
\address{Department of Mathematics\\ Ohio State University\\ 100 Math Tower, 231 West 18th Ave. \\ Columbus,  OH 43210\\ USA}
\email{hhtseng@math.ohio-state.edu}

\date{\today}

\begin{document}

\maketitle

\begin{abstract}
We derive a crepant resolution correspondence for some genus zero reduced Gromov-Witten invariants of Hilbert schemes of points on a K3 surface. 
\end{abstract}

\section{Introduction}
Let $S$ be a K3 surface. There are many references for K3 surfaces, e.g. \cite{Hu}.

\subsection{The geometry}
The symmetric group $S_n$ acts on the $n$-fold cartesian product $S^n$ by permuting the factors. The quotient
\begin{equation*}
S^n/S_n    
\end{equation*}
is an algebraic variety with quotient singularities. The ($n$-fold) symmetric product of $S$ is the stack quotient
\begin{equation*}
\mathsf{Sym}^n(S)=[S^n/S_n],    
\end{equation*}
which is a nonsingular projective Deligne-Mumford stack.

The Hilbert scheme of $n$ points on $S$,
\begin{equation*}
\mathsf{Hilb}^n(S),    
\end{equation*}
parametrizes $0$-dimensional closed subschemes of $S$ of length $n$. A standard reference for Hilbert schemes of points on surfaces is \cite{Na}. $\Hilb$ is smooth of dimension $2n$, see \cite[Theorem 1.15]{Na}.

The assignment $\mathsf{Hilb}^n(S)\ni [Z]\mapsto \sum_{x\in S} \text{length}(Z_x)[x]$ defines a morphism 
\begin{equation}
\pi: \mathsf{Hilb}^n(S)\to S^n/S_n,    
\end{equation}
see \cite[Theorem 1.5]{Na}. $\pi$ is an example of a Hilbert-Chow morphism. By \cite[Theorem 1.15]{Na}, $\pi$ is a resolution of singularities.

Since $S$ is a K3 surface, $S$ has a holomorphic symplectic form $\omega$. It follows tautologically that $\omega$ yields a holomorphic symplectic form on $\Sym$. By \cite[Theorem 1.17]{Na}, $\omega$ induces a holomorphic symplectic form on $\Hilb$.  It follows that both $\Hilb$ and $S^n/S_n$ have trivial canonical bundles. Thus $\pi$ is a {\em crepant} resolution of singularities\footnote{The Hilbert-Chow morphism is a crepant resolution for any smooth surfaces.}.

\subsection{Reduced Gromov-Witten theory}
Let $$\Mbar_{g,N}(\Hilb, \beta)$$ be the moduli stack of $N$-pointed genus $g$ stable maps to $\Hilb$ of class $\beta\in H_2(\Hilb, \mathbb{Z})\simeq H_2(S,\mathbb{Z})\oplus \mathbb{Z}$ (c.f. \cite[Corollary 3.13]{Ne1} and \cite[Section 0.2]{O}).  It carries a perfect obstruction theory with virtual dimension $(1-g)2n+3g-3+N$. Since the holomorphic symplectic form on $\Hilb$ is non-degenerate, the results of \cite{KL} implies that, when the $H_2(S,\mathbb{Z})$-factor of $\beta$ is nonzero, the virtual fundamental class vanishes, $$[\Mbar_{g,N}(\Hilb, \beta)]^{vir}=0.$$ The geometric reason behind this vanishing is the presence of a cosection of the obstruction sheaf of $\Mbar_{g,N}(\Hilb, \beta)$, which is constructed from the holomorphic symplectic form on $\Hilb$, see \cite{MP}. The cosection is surjective if the $H_2(S,\mathbb{Z})$-factor of $\beta$ is nonzero. By \cite{KL}, \cite{KL2}, \cite{O}, in this case $\Mbar_{g,N}(\Hilb, \beta)$ admits a {\em reduced} virtual fundamental class $$[\Mbar_{g,N}(\Hilb, \beta)]^{red}\in A_{(1-g)2n+3g-3+N+1}(\Mbar_{g,N}(\Hilb, \beta)).$$
The reduced Gromov-Witten theory of $\Hilb$ is defined by integrating against reduced virtual fundamental classes. For $\phi_1,...,\phi_N\in H^*(\Hilb)$, we define reduced Gromov-Witten invariants of $\Hilb$ to be 
$$\langle \prod_{i=1}^N\tau_{m_i}(\phi_i) \rangle_{g,\beta}^{\Hilb, red}:=\int_{[\Mbar_{g,N}(\Hilb, \beta)]^{red}}\prod_{i=1}^N\psi_i^{m_i}\ev_i^*(\phi_i).$$

Known calculations, e.g. \cite{O}, reveal very interesting structures present in the reduced Gromov-Witten theory of $\Hilb$.

Let $$\K_{g,N}(\Sym, \beta)$$ be the moduli stack of $N$-pointed genus $g$ orbifold stable maps to $\Sym$ with extended degree $\beta$ (in the sense of \cite[Section 2.1]{BG}). We assume that the  $H_2(S,\mathbb{Z})$-factor of $\beta$ is nonzero. The holomorphic symplectic form on $\Sym$ plays a similar role, implying that the virtual fundamental class vanishes, 
$$[\K_{g,N}(\Sym, \beta)]^{vir}=0.$$ 
Similarly, the obstruction sheaf of $\K_{g,N}(\Sym, \beta)$ admits a surjective cosection, leading to a reduced virtual fundamental class 
$$[\K_{g,N}(\Sym, \beta)]^{red}\in A_*(\K_{g,N}(\Sym, \beta)).$$
We define {\em reduced} Gromov-Witten theory\footnote{We can also define reduced Gromov-Witten invariants of $\Sym$ using twistor families \cite{BL}, \cite{R}.} of $\Sym$ by integrating against reduced virtual fundamental classes. For $\gamma_1,...,\gamma_N\in H^*({I}\Sym)$, we define reduced Gromov-Witten invariants of $\Sym$ to be $$\langle \prod_{i=1}^N\tau_{m_i}(\gamma_i) \rangle_{g,\beta}^{\Sym, red}:=\int_{[\K_{g,N}(\Sym, \beta)]^{red}}\prod_{i=1}^N\psi_i^{m_i}\ev_i^*(\gamma_i).$$

\subsection{Crepant resolutions}
The crepant resolution conjecture \cite{BG}, \cite{CIT}, \cite{CR}, concerns a correspondence between Gromov-Witten theory of a Gorenstein orbifold and Gromov-Witten theory of a crepant resolution of its underlying singular variety. For $(\mathsf{Hilb}^n(\mathbb{C}^2), \mathsf{Sym}^n(\mathbb{C}^2))$, crepant resolution conjecture is proven in full in \cite{PTs}, \cite{PTs2}. For $(\mathsf{Hilb}^n(Y), \mathsf{Sym}^n(Y))$ with $Y$ a toric del Pezzo surface, crepant resolution conjecture is proven for genus $0$ $3$-point invariants in \cite{Ne3}. For more general smooth projective surfaces $Y$, crepant resolution conjecture is proven in \cite{lq} for genus $0$ $3$-point invariants in {\em extremal} degrees, namely the sector of Gromov-Witten theory without curve classes from $Y$. 

For $(\Hilb, \Sym)$ with $S=K3$, crepant resolution conjecture is trivial because of the above vanishing. For these geometries, it is natural to consider crepant resolution conjecture for {\em reduced} Gromov-Witten theory.

Let $\mu=(\mu_1,...,\mu_k)$ be a partition of $n$ and let $\{\delta_1,...,\delta_{\text{rank}H^*(S)}\}$ be an additive basis of $H^*(S,\mathbb{Q})$. The collection 
\begin{equation*}
    \vec{\mu}=((\mu_1, \delta_{l_1}),...,(\mu_k, \delta_{l_k}))
\end{equation*}
is called a cohomology-weighted partition. Given $\vec{\mu}$, we associate a factor 
\begin{equation*}
\mathfrak{z}(\vec{\mu})=|\text{Aut}(\vec{\mu})|\prod_{i=1}^k \mu_i    
\end{equation*}
and a class $\delta_{l_1}\otimes...\otimes\delta_{l_k}\in H^*(S^\mu, \mathbb{Q})$. Here $S^\mu=S^{\ell(\mu)}$. Using the natural map $\bar{\pi}:\coprod_\mu S^\mu\to \bar{I}\Sym$ (see \cite[Section 2.5]{Ne3}), we define a class
\begin{equation*}
    \lambda(\vec{\mu}):=\frac{1}{\mathfrak{z}(\vec{\mu})}\bar{\pi}_*(\delta_{l_1}\otimes...\otimes\delta_{l_k})\in H^*(\bar{I}\Sym, \mathbb{Q}).
\end{equation*}
Given $\vec{\mu}$, we can also define a class in $H^*(\Hilb, \mathbb{Q})$ using Nakajima operators,
\begin{equation*}
    \theta(\vec{\mu}):=\frac{1}{\mathfrak{z}(\vec{\mu})}P_{\delta_{l_1}}[\mu_1]...P_{\delta_{l_k}}[\mu_k]\cdot 1\in H^*(\Hilb, \mathbb{Q}),
\end{equation*}
see \cite[Sections 2.7 and 6.2]{Ne3}, \cite[Chapter 8]{Na}, \cite[Section 0.2]{O} for details.

There is an additive isomorphism of graded vector spaces\footnote{Here $\text{age}(\mu)=n-\ell(\mu)$.}
\begin{equation}\label{eqn:coh_isom}
L: H^*({I}\Sym, \mathbb{Q}) \simeq H^*(\Hilb, \mathbb{Q}), \quad L(\lambda(\vec{\mu}))=(-\sqrt{-1})^{\text{age}(\mu)}\theta(\vec{\mu}),     
\end{equation}
see \cite[Proposition 6.2]{Ne3}, \cite[Proposition 3.5]{FG}, \cite{U}.

We consider the generating functions of genus $0$ invariants:
\begin{equation}\label{eqn:gen_func_hilb}
    \langle \gamma_1,\gamma_2,\gamma_3\rangle_{0,\gamma}^{\Hilb, red}(y):=\sum_{\mathsf{m}}\langle L(\gamma_1),L(\gamma_2),L(\gamma_3)\rangle_{0, (\gamma, \mathsf{m})}^{\Hilb, red}y^{\mathsf{m}},
\end{equation}
 
\begin{equation}\label{eqn:gen_func_sym}
    \langle \gamma_1,\gamma_2,\gamma_3\rangle_{0,\gamma}^{\Sym, red}(u):=u^{\sum_i \text{age}(\gamma_i)-2n}\sum_{\mathsf{h}}\langle \gamma_1,\gamma_2,\gamma_3\rangle_{0, (\gamma, \mathsf{h})}^{\Sym, red}u^{2\mathsf{h}-2}.
\end{equation}

The following is the main result of this paper.
\begin{thm}\label{thm:reduced_CRC}
Suppose $\gamma\in H_2(S,\mathbb{Z})$ is a nonzero class of divisibility at most $2$. Then the generating series $\langle \gamma_1,\gamma_2,\gamma_3\rangle_{0,\gamma}^{\Hilb, red}(y)$ is the Taylor expansion at $y=0$ of a rational function in $y$, and under the change of variables $-y=e^{\sqrt{-1} u}$, the following holds:
\begin{equation}\label{eqn:reduced_CRC}
\langle \gamma_1,\gamma_2,\gamma_3\rangle_{0,\gamma}^{\Hilb, red}(y)
= \langle \gamma_1,\gamma_2,\gamma_3\rangle_{0,\gamma}^{\Sym, red}(u).  
\end{equation}
\end{thm}

Theorem \ref{thm:reduced_CRC} is the crepant resolution correspondence for genus $0$ $3$-point primary reduced Gromov-Witten theories of $(\Hilb, \Sym)$. Following \cite{PTs2}, we also call this {\em $\mathsf{Hilb}/\mathsf{Sym}$ correspondence} for $S$.

\subsection{Calculations}
Calculations of reduced Gromov-Witten invariants of $\Hilb$ in \cite{O} can be combined with Theorem \ref{thm:reduced_CRC} to yield formulas for reduced Gromov-Witten invariants of $\Sym$. We discuss one example. 

Let $\pi: S\to \mathbb{P}^1$ be an elliptic K3 surface with a unique section $s: \mathbb{P}^1\to S$ and $24$ rational nodal fibers. In \cite[Section 2.2.1]{O}, $S$ is called a Bryan-Leung K3 surface.

Let $B$ be the class of a section of $\pi$ and $F$ the class of a fiber of $\pi$. For all $h\geq 0$, the class $\beta_h=B+hF\in H_2(S,\mathbb{Z})$ is a primitive and effective curve class. Furthermore, $\beta_h^2=2h-2$.

Consider cohomology-weighted partitions 
\begin{equation*}
\vec{\nu}:=(\underbrace{(1,F),...,(1,F)}_{n}) \quad \text{and}\quad  \vec{\eta}:=((1,F),\underbrace{(1,[pt]^\vee),...,(1,[pt]^\vee)}_{n-1}).
\end{equation*}
Here $(-)^\vee$ indicates Poincar\'e dual on $S$.

\begin{thm}\label{thm:cal_sym}
We have
\begin{equation*}
\sum_{h\geq 0}q^{h-1}\langle \lambda(\vec{\nu}), \lambda(\vec{\eta}), \lambda(\vec{\nu}) \rangle_{0, \beta_h}^{\Sym, red}(u)=\frac{(-\sqrt{-1})^{3n}}{(n!)^2}\frac{F(z,\tau)^{2n-2}}{\Delta(\tau)},    
\end{equation*}
where 
\begin{enumerate}
    \item $q=e^{2\pi\sqrt{-1}\tau}$ and $u=2\pi z$.

    \item $\Delta(\tau)=q\prod_{m\geq 1}(1-q^m)^{24}$ is the modular discriminant.

    \item $F(z,\tau)=(y^{1/2}+y^{-1/2})\prod_{m\geq 1}\frac{(1+yq^m)(1+y^{-1}q^m)}{(1-q^m)^2}$ is the Jacobi theta function ($y=-e^{2\pi\sqrt{-1}z}$).
\end{enumerate}
\end{thm}

\subsection{Outline}
Our proof of Theorem \ref{thm:reduced_CRC} is based on the work of D. Nesterov \cite{Ne1}, \cite{Ne2}, \cite{Ne3}. In \cite{Ne1}, \cite{Ne2}, Nesterov obtains, among other things, a correspondence between reduced Gromov-Witten theory of $\Hilb$ and reduced stable pair theory of a family of threefolds. In \cite{Ne3}, Nesterov obtains a correspondence between Gromov-Witten theory of a symmetric product and Gromov-Witten theory of a family of threefolds. These correspondences are obtained by means of wall-crossing formulas.

In Section \ref{subsec:review_Ne3} we review the results of \cite{Ne3}. In Section \ref{subsec:reduced_case}, we explain how to adjust \cite{Ne3} to obtain wall-crossing results in the reduced case. Theorem \ref{thm:reduced_CRC} is proved in Section \ref{subsec:pf_thm1} by combining the reduced wall-crossing results in \cite{Ne2} and Section \ref{subsec:reduced_case} and the Gromov-Witten/Pairs correspondence proven in \cite{O21}. Theorem \ref{thm:cal_sym} is proven in Section \ref{subsec:pf_thm2} by combining Theorem \ref{thm:reduced_CRC} and known calculations \cite{O}.

\begin{rem}[G. Oberdieck]
Since the Gromov-Witten/Pairs correspondence is proven for $\text{K3}\times C$ in curve classes of divisibility at most $2$ where $C$ is a smooth curve \cite[Theorem 1.2, Corollary 1.5]{O21}, the strategy for proving Theorem \ref{thm:reduced_CRC} should be applicable to derive a crepant resolution correspondence for $(\Hilb, \Sym)$ for genus one Gromov-Witten invariants with fixed $j$-invariant \cite[Equation (1)]{O21a}. Such a correspondence can then be combined with the calculation in \cite[Remark 1.2]{O21a} to obtain calculations of genus one Gromov-Witten invariants of $\Sym$ with fixed $j$-invariants.
\end{rem}

\subsection{Acknowledgment}
We thank G. Oberdieck and R. Pandharipande for helpful discussions. We are very grateful to the referee for important comments and suggestions. Throughout this work, D. G. is supported by OSU Speacial Graduate Assignment during his Ph.D. and later by AMS-Simons Travel Grant. H.-H. T. is supported in part by Simons Foundation Collaboration Grant. 

This work contains no data other than the paper itself. The authors are not aware of any conflict of interests.

\section{Gromov-Witten/Hurwitz wall-crossing}

\subsection{Review of Nesterov's work}\label{subsec:review_Ne3}
We summarize the main constructions of \cite{Ne3}. 

\subsubsection{Admissible maps}
Let $X$ be a smooth complex projective variety. Let $\epsilon\in (0,1]$. Let $$Adm_{g,N}^\epsilon(X^{(n)}, \beta)$$ 
be the stack of $\epsilon$-admissible maps as defined in \cite[Definition 2.5]{Ne3}. By \cite[Proposition 2.7]{Ne3} and \cite[Theorem 2.10]{Ne3}, $Adm_{g,N}^\epsilon(X^{(n)}, \beta)$ is a proper Deligne-Mumford stack of finite type. By \cite[Section 2.3]{Ne3}, $Adm_{g,N}^\epsilon(X^{(n)}, \beta)$ admits a perfect obstruction theory.

By the construction of \cite[Section 2.6]{Ne3}, there are evaluation maps  
taking values in a suitable variant of the inertia stack ${I}\mathsf{Sym}^n(X)$ of $\mathsf{Sym}^n(X)$. Descendant $\epsilon$-admissible invariants can be defined in the standard fashion, see \cite[Definition 2.12]{Ne3}.

\subsubsection{$\epsilon=1$}
When $\epsilon=1$, there is a map 
\begin{equation}\label{eqn:map_at_infty}
\rho: \K_{g,N}(\mathsf{Sym}^n(X), \beta)\to Adm_{g,N}^{1}(X^{(n)}, \beta),   
\end{equation}
see \cite[equation (9)]{Ne3}. By \cite[Lemma 2.11]{Ne3}, $\rho$ is virtually birational. By \cite[Lemma 2.13]{Ne3}, $1$-admissible invariants coincide with Gromov-Witten invariants of $\mathsf{Sym}^n(X)$.

\subsubsection{$\epsilon=0^+$}
When $\epsilon=0^+$, there is an identification \cite[equation (12)]{Ne3}
\begin{equation}\label{eqn:map_at_0}
Adm_{g,N}^0(X^{(n)}, \beta)=\Mbar_{\mathsf{m}}^\bullet(X\times C_{g,N}/\Mbar_{g,N}, (\gamma, n)),    
\end{equation}
where the right-hand side is a moduli space of relative stable maps with disconnected domains to the relative geometry
\begin{equation*}
   X\times C_{g,N}\to\Mbar_{g,N},  
\end{equation*}
where $C_{g,N}\to\Mbar_{g,N}$ is the universal family and relative divisors are given by marked points. This identification also identifies obstruction theories.

By \cite[Section 2.7.2]{Ne3}, $0^+$-admissible invariants coincide with relative Gromov-Witten invariants of $X\times C_{g,N}\to\Mbar_{g,N}$.

\subsubsection{Master space}\label{subsubsec:master_space}
The space $(0,1]$ is divided into chambers. When $\epsilon$ varies in a chamber, $Adm_{g,N}^\epsilon(X^{(n)}, \beta)$ stays the same. When $\epsilon$ crosses a wall between two chambers, $Adm_{g,N}^\epsilon(X^{(n)}, \beta)$ changes. The goal of wall-crossing is to study how $\epsilon$-admissible invariants change as $\epsilon$ crosses a wall. This is achieved in general in \cite[Section 4]{Ne3}, based on the geometry of master spaces introduced in \cite[Section 3]{Ne3}. 

Let $\epsilon_0$ be a wall. Let 
\begin{equation}\label{eqn:master_space}
    MAdm_{g,N}^{\epsilon_0}(X^{(n)}, \beta)
\end{equation}
be the moduli space of $N$-pointed genus $g$ $\epsilon_0$-admissible maps with calibrated tails, as defined in \cite[Definition 3.4]{Ne3}. By \cite[Theorem 3.5]{Ne3} and \cite[Theorem 3.9]{Ne3}, $MAdm_{g,N}^{\epsilon_0}(X^{(n)}, \beta)$ is a proper Deligne-Mumford stack of finite type. By \cite[Section 3.2]{Ne3}, $MAdm_{g,N}^{\epsilon_0}(X^{(n)}, \beta)$ admits a perfect obstruction theory, hence a virtual fundamental class.

The master space $MAdm_{g,N}^{\epsilon_0}(X^{(n)}, \beta)$ admits a $\mathbb{C}^*$-action obtained by scaling the $\mathbb{P}^1$ in the domain. The wall-crossing formula for $\epsilon$-admissible invariants \cite[Theorem 4.3]{Ne3} is obtained by applying virtual localization with respect to this $\mathbb{C}^*$-action, see \cite[Section 4.2]{Ne3} for more details.

\subsection{The reduced case}\label{subsec:reduced_case}
Now we consider the case $X=S$ is a K3 surface. We assume that the $H_2(S,\mathbb{Z})$-factor of $\beta$ is nonzero.

By the construction of \cite[Section 4.3--4.4]{MPT}, the holomorphic symplectic form on $S$ yields a cosection of the obstruction sheaf of $Adm_{g,N}^\epsilon(S^{(n)},\beta)$. Furthermore, the cosection is surjective. By \cite{KL}, \cite{KL2}, this yields a reduced virtual fundamental class $$[Adm_{g,N}^\epsilon(S^{(n)},\beta)]^{red}.$$ We can define reduced $\epsilon$-admissible invariants using the reduced virtual class:
\begin{equation}\label{def:reduced_adm_inv}
    \langle \prod_{i=1}^N\tau_{m_i}(\gamma_i) \rangle_{g,\beta}^{\epsilon, red}
\end{equation}
where $\gamma_1,...,\gamma_N\in H^*({I}\Sym)$.

When $\epsilon=1$, we have the reduced version of \cite[Lemma 2.11]{Ne3}:
\begin{equation}\label{eqn:reduced_vir_pushforward}
\rho_*[\K_{g,N}(\Sym, \beta)]^{red}=[Adm_{g,N}^{1}(S^{(n)},\beta)]^{red}.    
\end{equation}
Note that the reduced virtual classes are not constructed from perfect obstruction theories. Still by \cite{KL}, \cite{KL2}, reduced virtual classes are constructed by embedding intrinsic normal cones into the kernel of the cosection as cycles. Thus the proof of \cite[Theorem 5.0.1]{costello} (see also \cite{HW}) may be adapted to prove (\ref{eqn:reduced_vir_pushforward}). Therefore,
\begin{equation}\label{eqn:reduced_inv_sym}
 \langle \prod_{i=1}^N\tau_{m_i}(\gamma_i) \rangle_{g,\beta}^{1, red}=\langle \prod_{i=1}^N\tau_{m_i}(\gamma_i) \rangle_{g,\beta}^{\Sym, red}   
\end{equation}

Now we consider wall-crossing for reduced $\epsilon$-admissible invariants. As in Section \ref{subsubsec:master_space}, let $\epsilon_0\in (0,1]$ be a wall. Let $\epsilon_+, \epsilon_-$ be close to $\epsilon_0$ such that $\epsilon_+<\epsilon_0<\epsilon_-$. The master space $MAdm_{g,N}^{\epsilon_0}(S^{(n)}, \beta)$, as reviewed in Section \ref{subsubsec:master_space}, admits a surjective cosection (assuming the $H_2(S,\mathbb{Z})$-factor of $\beta$ is nonzero), again by the construction of \cite[Section 4.3--4.4]{MPT}. Hence it has a reduced virtual fundamental class $[MAdm_{g,N}^{\epsilon_0}(S^{(n)}, \beta)]^{red}$. By \cite{CKL}, we can apply virtual localization formula to the reduced virtual class. As in \cite[Section 4.2]{Ne3}, the $\mathbb{C}^*$-fixed loci of $MAdm_{g,N}^{\epsilon_0}(S^{(n)}, \beta)$ is described as follows:
\begin{equation}
MAdm_{g,N}^{\epsilon_0}(S^{(n)}, \beta)^{\mathbb{C}^*}=F_-\cup F_+\cup \coprod_{\vec{\beta}} F_{\vec{\beta}}.  
\end{equation}
As in \cite[Section 4.2.1]{Ne3}, we have $F_-=Adm_{g,N}^{\epsilon_-}(S^{(n)},\beta)$ and $[F_-]^{red}=[Adm_{g,N}^{\epsilon_-}(S^{(n)},\beta)]^{red}$. In the notation of \cite[Section 4.2.2]{Ne3}, we have $p_*[F_+]^{red}=[Adm_{g,N}^{\epsilon_+}(S^{(n)},\beta)]^{red}$. 

The fixed loci $F_{\vec{\beta}}$ are described as products, see \cite[Section 4.2.3]{Ne3} for details. We consider the reduced virtual class $[F_{\vec{\beta}}]^{red}$. Since the reduced virtual class of a product splits into a product of reduced and standard virtual classes (c.f. \cite[Section 3.9]{MPT}), and since standard virtual classes vanish whenever the cosection is surjective, which happens for maps representing nonzero curve classes, we see that $[F_{\vec{\beta}}]^{red}=0$ except for $\vec{\beta}=(\beta'=0,\beta)$. Thus we arrive at the reduced version of wall-crossing formula \cite[Theorem 4.3]{Ne3}:
\begin{equation}\label{eqn:reduced_wallcrossing}
\begin{split}
&\langle \prod_{i=1}^N\tau_{m_i}(\gamma_i) \rangle_{g,\beta}^{\epsilon_+, red}-\langle \prod_{i=1}^N\tau_{m_i}(\gamma_i) \rangle_{g,\beta}^{\epsilon_-, red}\\
=&\int_{[Adm_{g,N+1}^{\epsilon_-}(S^{(n)}, 0)]^{vir}}\prod_{i=1}^N\psi_i^{m_i}\ev_i^*(\gamma_i)\cdot \ev_{N+1}^*I_{+,\beta}(z)|_{z=-\psi_{N+1}}.    
\end{split}
\end{equation}
Here $I_{+,\beta}(z)$ is defined in \cite[Definition 4.1]{Ne3} with reduced virtual classes replacing the standard ones.

By a dimension argument similar to \cite[Section 5.1]{Ne3}, we see that $I_{+,\beta}(z)$ is proportional to the class $1\in H^0(\Sym)$. Since $Adm_{g,N+1}^{\epsilon_-}(S^{(n)}, 0)=\Mbar_{g,N+1}(\Sym,0)$, we see by string equation that $$\int_{[Adm_{0,N+1}^{\epsilon_-}(S^{(n)}, 0)]^{vir}}\prod_{i=1}^N\ev_i^*(\gamma_i)\cdot \ev_{N+1}^*I_{+,\beta}(z)|_{z=-\psi_{N+1}}=0,$$
whenever $2g-2+N>0$. Therefore, the reduced wall-crossing formula (\ref{eqn:reduced_wallcrossing}) without descendant insertions takes a simple form:
\begin{equation}\label{eqn:reduced_wallcrossing_primary}
\langle \prod_{i=1}^N\tau_{0}(\gamma_i) \rangle_{g,\beta}^{\epsilon_+, red}=\langle \prod_{i=1}^N\tau_{0}(\gamma_i) \rangle_{g,\beta}^{\epsilon_-, red}.   
\end{equation}

\section{Crepant resolution correspondence}

\subsection{Proof of Theorem \ref{thm:reduced_CRC}}\label{subsec:pf_thm1}
 Fix a nonzero $\gamma\in H_2(S,\mathbb{Z})$. Equations (\ref{eqn:reduced_wallcrossing_primary}) and (\ref{eqn:reduced_inv_sym}) imply that for $\gamma_1, \gamma_2, \gamma_3\in H^*({I}\Sym)$, 
\begin{equation}\label{eqn:reduced_GW_H}
\langle \gamma_1, \gamma_2, \gamma_3\rangle_{0,\beta}^{0^+, red}=\langle \gamma_1, \gamma_2, \gamma_3\rangle_{0,\beta}^{1, red}=\langle \gamma_1, \gamma_2, \gamma_3\rangle_{0,\beta}^{\Sym, red}.      
\end{equation}
Here, by construction, the left-hand side is the reduced relative Gromov-Witten invariants of $S\times \mathbb{P}^1$ relative to the divisors $S\times \{0, 1, \infty\}$.  

For $\Hilb$, \cite[Corollary 4.2]{Ne2} gives that for $\phi_1, \phi_2, \phi_3\in H^*(\Hilb)$,
\begin{equation}\label{eqn:reduced_Hilb_DT}
\langle \phi_1, \phi_2, \phi_3\rangle_{0,\beta}^{\Hilb, red}=\langle \phi_1, \phi_2, \phi_3\rangle_{n, \check{\beta}}^{S\times \mathbb{P}^1, red}.    
\end{equation}
Here the right-hand side is the reduced relative stable pair invariant of $S\times \mathbb{P}^1$ relative to the divisors $S\times \{0, 1, \infty\}$, see e.g. \cite{O21} for more details. 

The next step is to invoke proven Gromov-Witten/Pairs correspondence \cite{O21} for the relative geometry $(S\times \mathbb{P}^1, S\times \{0, 1, \infty\})$. 

Fix $\gamma_1,...,\gamma_N\in H^*({I}\Sym, \mathbb{Q})$ which are homogeneous with respect to the age grading.  We consider the following generating function of reduced $\epsilon$-quasimap invariants studied in \cite{Ne2}:
\begin{equation*}
    \langle \gamma_1,...,\gamma_N\rangle_{g,\gamma}^{\epsilon, red}(y):=\sum_{\mathsf{m}}\langle L(\gamma_1),...,L(\gamma_N)\rangle_{g, (\gamma, \mathsf{m})}^{\epsilon, red}y^{\mathsf{m}}.
\end{equation*}
We also consider the following generating series of reduced $\epsilon$-admissible invariants:
\begin{equation*}
  \langle \gamma_1,...,\gamma_N\rangle_{g,\gamma}^{\epsilon, red}(u):=u^{\sum_i \text{age}(\gamma_i)-2n}\sum_{\mathsf{h}}\langle \gamma_1,...,\gamma_N\rangle_{g, (\gamma, \mathsf{h})}^{\epsilon, red}u^{2\mathsf{h}-2}.
\end{equation*}
In these notations, the reduced Gromov-Witten/Pairs correspondence for the family 
\begin{equation}\label{eqn:family_3fold}
S\times C_{g,N}\to \Mbar_{g,N}
\end{equation}
 can be formulated as follows, c.f. \cite[Section 6.3]{Ne3}: 
$\langle \gamma_1,...,\gamma_N\rangle_{g,\gamma}^{0^+, red}(y)$ is the Taylor expansion at $y=0$ of a rational function, and under the change of variables $-y=e^{\sqrt{-1} u}$, the following holds
\begin{equation}\label{eqn:GW_P_general}
 \langle \gamma_1,...,\gamma_N\rangle_{g,\gamma}^{0^+, red}(y)
 =\langle \gamma_1,...,\gamma_N\rangle_{g,\gamma}^{0^+, red}(u).   
\end{equation}

The correspondence (\ref{eqn:GW_P_general}) is still conjectural in general. When $g=0, N=3$, the family (\ref{eqn:family_3fold}) reduces to $S\times \mathbb{P}^1$. By \cite[Theorem 1.2, Corollary 1.5]{O21}, (\ref{eqn:GW_P_general}) holds when $g=0$, $N=3$, and $\gamma\in H_2(S,\mathbb{Z})$ is of divisibility at most $2$. Together with (\ref{eqn:reduced_GW_H}) and (\ref{eqn:reduced_Hilb_DT}), we obtain Theorem \ref{thm:reduced_CRC}.

\begin{rmk}
To prove the extension of Theorem \ref{thm:reduced_CRC} when $(g,N)\neq (0,3)$, we need to prove reduced Gromov-Witten/Pairs correspondences for families (\ref{eqn:family_3fold}) when the base $\Mbar_{g,N}$ is not a point. To the best of our knowledge, such a correspondence is not known for K3 surfaces $S$: in fact, to date, Gromov-Witten/Pair correspondences for (\ref{eqn:family_3fold}) is only known for surface $\mathbb{C}^2$ and is established by proving $\mathsf{Hilb}/\mathsf{Sym}$ correspondence for $\mathbb{C}^2$, see \cite{PTs}.   
\end{rmk}

\subsection{Proof of Theorem \ref{thm:cal_sym}}\label{subsec:pf_thm2} 
The generating function
\begin{equation}
 \langle \theta(\vec{\nu}), \theta(\vec{\nu}) \rangle^{\Hilb}(y,q)=\sum_{h\geq 0}\sum_{k\in \mathbb{Z}}y^kq^{h-1}\langle \theta(\vec{\nu}), \theta(\vec{\nu}) \rangle_{0, (\beta_h, k)}^{\Hilb, red}
\end{equation}
is evaluated in \cite[Theorem 10]{O} as follows:
\begin{equation}\label{eqn:cal_hilb}
\langle \theta(\vec{\nu}), \theta(\vec{\nu}) \rangle^{\Hilb}(y,q)=\frac{F(z,\tau)^{2n-2}}{(n!)^2\Delta(\tau)},    
\end{equation}
where $q=e^{2\pi\sqrt{-1}\tau}$, $y=-e^{2\pi\sqrt{-1}z}$. Here the factor $(1/n!)^2$ appears because $\theta(\vec{\nu})=\frac{1}{n!}P_F[-1]^n\cdot 1$.

Since we have 
\begin{equation*}
\int_{\beta_h+kA} \theta(\vec{\eta})=1,   
\end{equation*}
divisor equation implies, 
\begin{equation}
\sum_{k\in \mathbb{Z}}y^k\langle \theta(\vec{\nu}), \theta(\vec{\eta}), \theta(\vec{\nu}) \rangle_{0, (\beta_h, k)}^{\Hilb, red}=\sum_{k\in \mathbb{Z}}y^k\langle \theta(\vec{\nu}), \theta(\vec{\nu}) \rangle_{0, (\beta_h, k)}^{\Hilb, red}.
\end{equation}
By (\ref{eqn:reduced_CRC}), we have 
\begin{equation}
\begin{split}
&\sum_{k\in \mathbb{Z}}y^k\langle \theta(\vec{\nu}), \theta(\vec{\eta}), \theta(\vec{\nu}) \rangle_{0, (\beta_h, k)}^{\Hilb, red}\\
=&\langle L^{-1}\theta(\vec{\nu}), L^{-1}\theta(\vec{\eta}), L^{-1}\theta(\vec{\nu}) \rangle_{0, \beta_h}^{\Sym, red}(u)\\
=&\frac{1}{(-\sqrt{-1})^{3n}}\langle \lambda(\vec{\nu}), \lambda(\vec{\eta}), \lambda(\vec{\nu}) \rangle_{0, \beta_h}^{\Sym, red}(u),
\end{split}
\end{equation}
after $-y=e^{\sqrt{-1}u}$. Theorem \ref{thm:cal_sym} follows by combining this with (\ref{eqn:cal_hilb}).

\end{document}